\documentstyle[twoside, 12pt, amstex]{amsart}

\makeatletter

\begin{document}

\newcommand{\e}{\epsilon}
\newcommand{\w}{{\bold w}}
\newcommand{\y}{{\bold y}}
\newcommand{\z}{{\bold z}}
\newcommand{\x}{{\bold x}}
\newcommand{\N}{{\mathbb N}}
\newcommand{\Z}{{\mathbb Z}}
\newcommand{\F}{{\bold F}}
\newcommand{\R}{{\mathbb R}}
\newcommand{\Q}{{\mathbb Q}}
\newcommand{\Ql}{{\mathbb Q_{\ell}}}
\newcommand{\C}{{\mathbb C}}
\newcommand{\K}{{\mathbb K}}
\newcommand{\BP}{{\mathbb P}}
\newcommand{\cF}{{\mathcal F}}
\newcommand{\cL}{{\mathcal L}}
\newcommand{\cR}{{\mathcal R}}
\newcommand{\cO}{{\mathcal O}}
\newcommand{\cX}{{\mathcal X}}
\newcommand{\cH}{{\mathcal H}}
\newcommand{\cM}{{\mathcal M}}
\newcommand{\sB}{{\sf B}}
\newcommand{\cT}{{\mathcal T}}
\newcommand{\cI}{{\mathcal I}}
\newcommand{\cS}{{\mathcal S}}
\newcommand{\sE}{{\sf E}}

\newcommand{\sA}{{\sf A}}
\newcommand{\ga}{{\sf a}}
\newcommand{\es}{{\sf s}}
\newcommand{\m}{{\bold m}}
\newcommand{\bS}{{\bold S}}
\newcommand{\Fp}{{\mathbb F}_p}
\newcommand{\Fq}{{\mathbb F}_q}

\newcommand{\ihra}{\stackrel{i}{\hookrightarrow}}
\newcommand\rank{\mathop{\rm rank}\nolimits}
\newcommand\im{\mathop{\rm Im}\nolimits}
\newcommand\Li{\mathop{\rm Li}\nolimits}
\newcommand\NS{\mathop{\rm NS}\nolimits}
\newcommand\Hom{\mathop{\rm Hom}\nolimits}
\newcommand\Pic{\mathop{\rm Pic}\nolimits}
\newcommand\Spec{\mathop{\rm Spec}\nolimits}
\newcommand\Hilb{\mathop{\rm Hilb}\nolimits}
\newcommand{\length}{\mathop{\rm length}\nolimits}

\newcommand\lra{\longrightarrow}
\newcommand\ra{\rightarrow}
\newcommand\cJ{{\mathcal J}}
\newcommand\JG{J_{\Gamma}}
\newcommand{\wvskp}{\vspace{1cm}}
\newcommand{\vskp}{\vspace{5mm}}
\newcommand{\nvskp}{\vspace{1mm}}
\newcommand{\nid}{\noindent}
\newcommand{\new}{\nvskp \nid}
\newtheorem{Assumption}{Assumption}[section]
\newtheorem{Theorem}{Theorem}[section]
\newtheorem{Lemma}{Lemma}[section]
\newtheorem{Remark}{Remark}[section]
\newtheorem{Corollary}{Corollary}[section]
\newtheorem{Conjecture}{Conjecture}[section]
\newtheorem{Proposition}{Proposition}[section]
\newtheorem{Example}{Example}[section]
\newtheorem{Definition}{Definition}[section]
\newtheorem{Question}{Question}[section]
\newtheorem{Problem}{Problem}[section]

\renewcommand{\thesubsection}{\it}

\title{ The modularity conjecture for rigid Calabi--Yau threefolds
over $\Q$} 

\author{Masa-Hiko Saito and Noriko Yui}
\thanks{M.-H. Saito was partially supported by JSPS and Grant-in Aid
for Scientific Research (B-09440015), (B-12440008), 
the Ministry of
Education, Science and Culture, Japan.  N. Yui was partially
supported by a Research Grant from NSERC, Canada.}
\address{Department of Mathematics, Faculty of Science, 
Kobe University, Rokko, Kobe, 657-8501, Japan}
\email{mhsaito@math.kobe-u.ac.jp}
\address{Department of Mathematics, Queen's University, Kingston.
Ontario         Canada K7L 3N6}

\email{ yui@mast.queensu.ca;
yui@fields.utoronto.ca}
\date{September 14, 2000}

\begin{abstract}: We formulate the modularity conjecture for rigid Calabi--Yau 
threefolds defined over the field $\Q$ of rational numbers.  
We establish the modularity for the rigid Calabi--Yau threefold arising 
from the root lattice $A_3$. Our proof is based on geometric analysis. 
\end{abstract}

\maketitle

\section{\bf The $L$--series of Calabi--Yau 
threefolds}

 Let $\Q$ be the field of rational numbers, and let $\bar{\Q}$ be its 
algebraic closure with Galois group ${\mathcal{G}}:=\text{Gal}(\bar{\Q}/\Q)$.  
Let $X$ be a smooth projective threefold defined over $\Q$ or more generally
over a number field.

\begin{Definition}{\rm 
 $X$ is a {\it Calabi--Yau} threefold
if it satisfies the following two conditions: 

 (a) $H^1(X,\cO_X)=H^2(X,\cO_X)=0$, and 

 (b) The canonical bundle is trivial, i.e., $K_X\simeq\cO_X$.
}
\end{Definition}

\vspace{0.2cm}
\noindent{\bf The numerical invariants of Calabi--Yau
threefolds.}

\vspace{0.2cm}
Let $X$ be a Calabi--Yau threefold defined 
over $\Q$, and let
$\bar X =X \times_{\Q} \bar {\Q}$. 

\vspace{0.2cm}
(1) The $(i,j)$-th {\it Hodge} number $h^{i,j}(X)$ of $X$ is defined
by 
$$h^{i,j}(X)=\text{dim}_{\bar \Q} H^j(\bar X, \Omega_{\bar X}^i).$$
There are symmetries among Hodge numbers: 
$h^{i,j}(X)=h^{j,i}(X)$, and $h^{i,j}(X)=h^{3-j,3-i}(X).$ 

The first identity follows from the complex conjugation operation
and the latter from Serre duality on the Hodge cohomology groups.  
The condition (a) implies that  $h^{1,0}(X)=h^{2,0}(X)=0$, and 
the condition (b) that $h^{3,0}(X)=h^{0,3}(X)=1$.
The number $h^{2,1}(X)$ represents the number of deformations of 
complex structures on $X$, and $h^{1,1}(X)$ is the number of Hodge 
$(1,1)$-cycles on $X$.  The Hodge numbers are concocted to define 
the {\it Hodge diamond}:
 
$$\matrix 
1\\
0\quad\quad \quad 0\\
0 \quad\quad  \,\,h^{1,1}\quad\quad \quad 0\\
1\quad\quad\quad h^{1,2} \quad\quad h^{1,2} \quad\quad\quad 1\\
0 \quad\quad \,\,h^{2,2}\quad\quad\quad 0 \\
0\quad\quad\quad \quad 0\\
1
\endmatrix
$$

\vspace{0.5cm}
(2) The $n$-th {\it Betti} number of $X$  is defined by
$$B_n(X)=\text{dim}_{\Ql} H^n_{\text{et}}(\bar X,\Ql)=\text{dim}_{\C}
H^n(X\otimes\C, \C)$$
where the cohomology groups are, 
respectively,
the $\ell$-adic \'etale and singular cohomology groups.   
The Hodge decomposition $H^n(X \otimes \C, \C) = 
\oplus_{i+j=n} H^j(X \otimes \C, \Omega^i_{X \otimes \C})$ ensures that  
$$B_n(X)=\sum_{i+j=n} h^{i,j}(X)=\sum_{i=0}^n h^{i,n-i}(X).$$
Using the above formula and the Poincar\'e duality $B_n(X)=B_{2 \text{dim}(X)-n}$, 
we can compute 
the Betti numbers of Calabi--Yau threefolds:
$$B_0(X)=1=B_6(X),\quad B_1(X)=0=B_5(X),$$
$$B_2(X)=h^{1,1}(X)=h^{2,2}(X)=B_4(X)\quad\text{and}\quad B_3(X)=2(1+h^{2,1}(X)).$$
The (topological) Euler characteristic is $\chi(X)=2(h^{1,1}(X) - h^{2,1}(X)).$
\medskip

\begin{Definition} 
{\rm A smooth projective Calabi--Yau threefold $X$ over $\Q$ (or over any 
field) is called {\it rigid} if $h^{2,1}(X)=h^{1,2}(X)=0$ so that $B_3(X)=2$.}
\end{Definition}

\begin{Definition}{\rm  Let $X$ be a Calabi--Yau threefold
defined over $\Q$. Assume that $X$ has a suitable integral model.
The $L$--series of $X$ is defined to be the $L$--series
of the (semi-simplification of the) Galois representation on
$H^3_{\text{et}}(\bar{X},\Ql)$.  That is,
$$
L(X,s):=L(H^3_{\text{et}}(\bar{X},\Ql),s).
$$}
\end{Definition}

\vspace{0.5cm}
\noindent{\bf Digression:} Let $X$ be a Calabi--Yau threefold defined 
over $\Q$ with a suitable integral model.
Let $p$ be a {\it good} prime for $X$, that is, the reduction $X\pmod p$
defines a smooth projective variety over the prime field $\Fp$.
Let $\text{Frob}_p$ denote the (geometric) Frobenius morphism of $X$ at $p$. 
We consider the action of $\text{Frob}_p$ on the $\ell$-adic etale cohomology
group $H^3_{\text{et}}(\bar{X},\Ql)$, and let $t_3(p)$ denote the
trace of the $\text{Frob}_p$ on $H^3_{\text{et}}(\bar{X},\Ql)$. 
By the Lefschetz fixed point formula, $t_3(p)$ can be determined by
counting the number of $\Fp$--rational points on $X$:
$$t_3(p)=1+p^3+(1+p)t_2(p)-\#X(\Fp).$$
Here $t_2(p)=p\,\text{Trace}(\rho(\text{Frob}_p))$ where
$\rho$ is an Artin representation of dimension $h^{1,1}(X)$  so in
particular its trace is an integer which is a sum of
$h^{1,1}(X)$ roots of unity.  This implies that
$|t_2(p)|\leq p\,h^{1,1}(X)$, and for a set of primes of
positive density, we get $t_2(p)=p h^{1,1}(X)$. 
Define
$$P_{3,p}(T):=\text{det}(1-\text{Frob}_p\,T\,|\, H^3_{\text{et}}(\bar{X},\Ql)).$$

If $X$ is a rigid Calabi--Yau threefold, then $P_{p,3}(T)$ is an integral 
polynomial of degree $\text{deg}(P_{3,p})=2$; it is of the form:
$$P_{3,p}(T)=1-t_3(p)T+p^3T^2\in 1+T\Z[T]$$
where  $t_3(p)$ is subject to the Riemann Hypothesis: $|t_3(p)|\leq 2p^{3/2}$.  
The $L$--series \newline $L(H^3_{\text{et}}(\bar{X},\Ql), s)$
is then given by
$$L(H^3_{\text{et}}(\bar{X}, \Ql),s)=(*)\prod_{p} P_{3,p}(p^{-s})^{-1}$$
where $p$ runs over good primes, and $(*)$ is the Euler factor corresponding to
bad primes.
\medskip

\section{\bf {The modularity conjecture for 
rigid Calabi--Yau threefolds over $\Q$}}

 Let $ k \geq 1$ be an integer.  Let $\Gamma$ be an arithmetic subgroup of 
 $\text{SL}_2(\Z)$ of finite index.  We denote by $S_k(\Gamma)$ the complex
vector space of all cusp forms of weight $k$ with respect to $\Gamma$.
We now formulate the modularity conjecture for rigid Calabi--Yau
threefolds defined over $\Q$.
\medskip
 
\begin{Conjecture}\label{mconjecture}{\bf The modularity conjecture}: 
{\rm Any rigid Calabi--Yau threefold $X$ defined over $\Q$ is modular in
the sense that its $L$--series of $X$ coincides with the Mellin transform of the 
$L$--series of a modular (cusp) form $f$ of weight $4$ on $\Gamma_0(N)$. 
Here $N$ is a positive integer divisible by the primes of bad reduction.  
More precisely, we have, up to a finite Euler factors, 
\begin{equation*}
L(X,s)=L(f,s)\quad\text{for}\quad f\in S_4(\Gamma_0(N)).
\end{equation*}} 
\end{Conjecture}
\smallskip

Here are some justifications for formulating the modularity conjecture
for rigid Calabi--Yau threefolds over $\Q$.
\smallskip

\begin{Remark}{\rm The conjecture of Taniyama--Shimura--Weil that every elliptic 
curve over $\Q$ is modular, has been established by Wiles and his former students
in totality (see Breuil, Conrad, Diamond and Taylor \cite{BCDT}).
Noting that elliptic curves are dimension one Calabi--Yau varieties,
our modularity conjecture (2.1) may be regarded a dimension three
generalization of the Taniyama--Shimura--Weil conjecture to rigid Calabi--Yau
threefolds over $\Q$.}
\end{Remark} 
\smallskip

\begin{Remark}{\rm  Livn\'e \cite{Li2} considered a rank $2$ motive 
$\mathcal{M}$ over $\Q$ with Hodge numbers 
$h^{p,q}({\mathcal{M}})=h^{q,p}({\mathcal{M}})=1$ 
where $p>q$, which respects an orthogonal form up to similitudes.  
Livn\`e showed how to express the $L$--series of ${\mathcal{M}}$ 
in terms of Hecke characters.

The examples of rank $2$ motives ${\mathcal{M}}$ considered by Livn\`e arose 
from an elliptic curve $E$ over $\Q$
or from a singular $K3$ surface $X$, i.e., ${\mathcal{M}}=H^1(E)\otimes \Q$
or $H^2(X,\Z)/\text{Pic}(X)\otimes \Q$.  In the latter case,
Shioda and Inose \cite{SI} determined the $L$--series of a singular 
$K3$ surface $X$ (up to a finite numbers of 
Euler factors) passing to some finite extension of $\Q$: 
\begin{equation*}
L(X,s) = L(\psi_1^2,s)L(\bar\psi_1^2,s)
\end{equation*} 
where $\psi_1$ is the Hecke character associated to an elliptic 
curve with complex multiplication.} 
\end{Remark}
\smallskip

\begin{Remark}{\rm  Fontaine and Mazur \cite{FM} have conjectured 
that all irreducible
odd $2$--dimensional Galois representations ``coming from geometry''
should be modular, up to a Tate twist.  Our modularity conjecture
(2.1) may be regarded as a concrete realization of the Fontaine--Mazur
conjecture.  For a recent progress on a conjecture of Fontaine--Mazur,
the reader is referred to a paper of Richard Taylor \cite {T}.}
\end{Remark}
\smallskip

\begin{Remark} {\rm Let $X$ be a projective variety of odd dimension $m$ 
over $\Bbb Q$ such that the $m$-th Betti number
$B_m={\rm{dim}}\, H^m(X\otimes{\Bbb C},{\Bbb C})=2$ and 
$H^m(X\otimes{\Bbb C},{\Bbb C})$ has the Hodge decomposition of 
type $(m,0)+(0,m)$. Serre \cite{Sr} has formulated a modularity conjecture for the
residual mod $p$ $2$-dimensional Galois representation attached to $X$ for all
primes.  In particular, our modularity conjecture 2.1 is 
a special case of the conjectural Theorem 6 of Serre \cite{Sr} for $m=3$.
Serre has informed us about this in his email dated
June 23, 2000.  We are thankful to him for pointing this out.}
\end{Remark}
\medskip

\section{\bf A rigid Calabi--Yau threefold arising from the
root lattice $A_3$}

\hskip 1cm This example of a rigid  Calabi--Yau threefold is 
constructed from the root lattice $A_3$ via a toric construction.

\vspace{0.5cm}
{\bf {(3.1) Toric construction}.}
For the general backgrounds on toric varieties, the reader is
referred to Batyrev \cite{Bat}, Fulton \cite{Fu} and 
Fulton and Harris \cite{FH}.

We consider a root system $\cR$ of rank $r$.  Let $\cL_{\cR}$
be the root lattice generated by $\cR$, and let $\cL_{\cR}^*$ be its 
dual lattice. Let $\Sigma_{\cR}$ be the fan in $\cL_{\cR}^*\otimes\Q$.
The fan $\Sigma_{\cR}$ gives rise to a toric variety, which is denoted 
by $X(\Sigma_{\cR})$.  Let $\Delta_{\cR}$ be a polyhedron in $\cL_{\cR}\otimes\Q$ 
with vertices in $\cR$, and let $L(\Delta_{\cR})$ denote the space of Laurent
polynomials with support in $\Delta_{\cR}$.  For each root
$r\in\cR$, assign a monomial ${\bold e}^r$, and define a Laurent polynomial
$$\chi_{\cR}:=\sum_{r \in \cR} {\bold e}^r \in L(\Delta_{\cR}).$$
Then we have a rational function $\chi_{\cR}:X(\Sigma_{\cR})\to \BP^1$.
For each $\lambda \in \BP^1$, define
$$X_\lambda:=\overline{\chi_{\cR}^{-1}(\lambda)}\subset X(\Sigma_{\cR}).$$\
There is a base locus $\mathcal B$.  Let
$\mathcal B:X_{\cR}\to X(\Sigma_{\cR})$
be the blow up of the base locus.  Then we obtain a pencil of 
varieties
$$\cX_{\cR}:=\{(\lambda,x)\in \BP^1\times X(\Sigma_{\cR})\,|\,
\overline{\chi_{\cR}(x)=\lambda}\}.$$
We are interested in root lattices $\cR$ which give
rise to Calabi--Yau threefolds $\cX_{\cR}$.
\medskip

{\bf {(3.2) The root lattice $A_3$ and the construction of
a Calabi--Yau threefold}.}
Verrill\cite{V} constructed a Calabi--Yau threefold associated
to the root lattice $A_3$.  We will briefly describe
Verrill's construction.  Let $\{E_i\,|\,i=1,2,3,4\}$ be the standard
basis for $\R^4$.  The root lattice $A_3$ is a sublattice
of $\R^4$ of rank $3$ generated by
 $v_1:=E_1-E_2,\,v_2:=E_2-E_3,\,v_3:=E_3-E_4$
and the collection of all roots is given by the set
$$\{E_i-E_j\,|\,1\leq i,j\leq 4,\, i\neq j\}.$$
To a root $E_i-E_j$, we associate a monomial $X_iX_j^{-1}$ by putting
$X_i={\bold e}^{E_i}$.
Then the character of the adjoint representation is given by
$$\chi_{A_3}=X_1X_2^{-1}+X_1X_3^{-1}+X_1X_4^{-1}+X_2X_3^{-1}+X_2X_4^{-1}+
X_3X_4^{-1}$$
$$+X_2X_1^{-1}+X_3X_1^{-1}+X_4X_1^{-1}+X_3X_2^{-1}+X_4X_2^{-1}+X_4X_3^{-1}$$
$$=(X_1+X_2+X_3+X_4)(X_1^{-1}+X_2^{-1}+X_3^{-1}+X_4^{-1})-4.$$ 
Let $X(\Sigma_{A_3})\to \BP^1$ denote the
toric variety defined by $\chi_{A_3}=\lambda \in \BP^1$. It is a rational variety.  
Now we take a double cover  of $\cX(\Sigma_{A_3})$ by
putting $\lambda=\frac{(t-1)^2}{t}$ with $t\in \BP^1$. 
Then a desingularization $X$ of the double cover is a Calabi--Yau 
threefold defined over $\Q$ whose generic fiber $X_{\lambda}$ is a $K3$ surface.
\medskip

\begin{Lemma} The numerical characters of $X$ are given 
as follows:
$$h^{1,2}=0,\,h^{3,0}=1,\,h^{1,1}=50,\, \chi(X)=100$$
and hence $X$ is a rigid Calabi--Yau threefold.
\end{Lemma}
\smallskip

We can see from the defining equation of $X$ that $2$ and $3$
are bad primes.
For this rigid Calabi--Yau threefold $X$ over $\Q$, we can establish the 
modularity conjecture.  We will give a proof of the following theorem 
in the section 6 below. Verrill \cite{V} also proved the theorem with 
totally different method (\`a la Livn\'e \cite{Li1}) from ours. 
\smallskip

\begin{Theorem}\label{thm:modular} Let $X$ be the rigid
Calabi--Yau threefold constructed by Verrill.  Then
$L$--series of $X$ coincides, up to the Euler factors corresponding
to primes $2$ and $3$, with the Mellin
transform of a weight $4$ newform on $\Gamma_0(6)$. That is,
$$
L(X,s) =
L(f,s) \quad\text{for}\quad f\in S_4(\Gamma_0(6)).$$
The cusp form $f$ has the $q$--expansion of the form:
$$
f(q)=\eta(q)^2\eta(q^2)^2\eta(q^3)^2\eta(q^6)^2.
$$ 
\end{Theorem}
\medskip

\section{\bf The self-products of elliptic modular 
surfaces.}

\hskip 1cm Let $\Gamma$ be an arithmetic subgroup of $\text{SL}_2({\Z})$ 
with no torsion elements and let $\cH$ be the upper half complex plane.  
Then we can define the modular curve $C_{\Gamma}^0 :=\cH/\Gamma$
and its compactification $C_{\Gamma}:=\overline{\cH/\Gamma}
=\cH/\Gamma\cup\{\rm{cusps\,\, of}\,\,\Gamma\}$. 
It is known (cf. \cite{D}) that
the universal family of elliptic curves 
$$f^0: S^{0}_{\Gamma} \lra C_{\Gamma}^0 = \cH/\Gamma$$
exists, and that its N\'eron model has the unique 
smooth minimal compactification
$$ f:= f_{\Gamma}:S_{\Gamma} \lra C_{\Gamma}= \overline{\cH/\Gamma}.$$  
The fibration  $f: S_{ \Gamma } \lra C_{\Gamma}$
is called the {\it elliptic modular surface associated to $\Gamma$.}   
\smallskip

Let us now recall some basic facts about elliptic modular surfaces,
which are relevant to our subsequent discussions. 
Consider the following diagram
$$
\begin{array}{ccc}
S_{\Gamma}^0 & \hookrightarrow & S_{\Gamma} \\
f^0 \downarrow \quad  &     &  f \downarrow \quad \\
C_{\Gamma}^0  & \stackrel{j}{\hookrightarrow} & C_{\Gamma}. 
\end{array}
$$
Making use of the invariant cycle theorem, 
we obtain an isomorphism of sheaves
$$
R^1 f_* \Q_{S_{\Gamma}}
\simeq 
j_{*} R^1 f^0_* \Q_{S_{\Gamma}^0}.   
$$
In fact, we have the natural map 
\begin{equation}\label{eq:inv}
R^1 f_* \Q_{S_{\Gamma}}
\lra 
j_{*} R^1 f^0_* \Q_{S_{\Gamma}^0}
\end{equation}
which is an isomorphim over 
$ C_{\Gamma}^0 = C_{\Gamma}\setminus \Sigma $  
where $\Sigma = \{p_1, \cdots, p_r \} $ is the set of cusps of $\Gamma$. 
Looking at the stalks at a cusp $p_i$, the morphism 
\begin{equation}\label{eq:stalk}
(R^1 f_* \Q_{S_{\Gamma}})_{p_i}
\lra 
(j_{*} R^1 f^0_* \Q_{S_{\Gamma}^0})_{p_i}
\end{equation}
is surjective by local invariant cycle theorem (cf. e.g. [Proposition 15.12, \cite{Z}]). 
On the other hand, it is easy to see that both of the stalks in (\ref{eq:stalk})  are 
isomorphic to a one-dimensional $\Q$-vector space.  
Hence the linear map (\ref{eq:stalk}) is an isomorphism.  

\smallskip

Let $V$ denote the natural representation of $SL_2(\Q)$, and let 
$S^k(V)$ be the $k$-th symmetric tensor representation of $V$.   
Then it is known (\cite{Z}, Proposition 12.5) that for each
$k\geq 0$ there exists 
an isomorphism 
\begin{equation}\label{eq:paraiso}
\tilde{H}^1(\Gamma, S^k(V)) \simeq 
H^1(C_{\Gamma}, j_*(S^k(R^1 f^0_*\Q_{S_{\Gamma}^0}))) \quad\rm{for\,\, each}
\quad k\geq 0.
\end{equation}
Here $\tilde{H}^1(\Gamma, S^k(V))$ is the 
parabolic cohomology group associated to the 
representation of $\Gamma$ on $S^k(V)$. 
Let $\omega_{S_{\Gamma}/C_{\Gamma}} = K_{S_{\Gamma}} 
\otimes f^{*}(K_{C_{\Gamma}}^{-1})$ be the relative canonical bundle of $f$,
and we set  
$$\omega=f_*(\omega_{S_{\Gamma}/C_{\Gamma}}).
$$
Then we see that $\omega$ is an invertible sheaf on $C_{\Gamma}$.  
\smallskip

Now we are ready to describe an isomorphism, the so-called {\it Shimura
isomorphism} and its far-reaching consequences.  

\begin{Proposition} \rm{(\cite{Sh}, \cite{D}, \cite{Z}.)}
{\sl There exists a natural isomorphism (\it {the Shimura isomorphism})
\begin{equation}\label{eq:Shimura}
H^0(C_{\Gamma}, \Omega^1_{C_{\Gamma}} \otimes 
\omega^k) \simeq S_{k+2}(\Gamma)
\end{equation}
and it gives rise to the following commutative diagram:
$$
\begin{array}{ccc}
\tilde{H}^1(\Gamma, S^k(V))  \otimes_{\Q} \C &  \simeq & 
  S_{k+2}(\Gamma) \oplus \overline{S_{k+2}(\Gamma)} \\
\downarrow  &   &  \downarrow \\
 H^1(C_{\Gamma},  j_*(S^k(R^1 f^0_*\C_{S_{\Gamma}^0}))) &  \simeq &  
 H^0(C_{\Gamma}, \Omega^1_{C_{\Gamma}} \otimes \omega^k) \oplus  
 \overline{H^0(C_{\Gamma}, \Omega^1_{C_{\Gamma}} \otimes \omega^k)}
\end{array}
$$
The first vertical arrow is the isomorphism 
induced by (\ref{eq:paraiso}).} 
\end{Proposition}
\smallskip

With the Shimura isomorphism at our disposal, we
may now characterize rational elliptic modular surfaces.

\begin{Proposition}
The elliptic modular surface associated to $\Gamma$, $f:S_{\Gamma} \lra C_{\Gamma}$ 
is a rational elliptic surface, if and only if
$$\dim S_2(\Gamma) = \dim S_3(\Gamma) = 0.$$  
\end{Proposition}
\smallskip

{\it Proof.} From the formula 
for the canonical bundle of the elliptic fibration
$f:S_{\Gamma} \lra C_{\Gamma}$, we see that 
\begin{equation}\label{eq:canonical}
K_{S_{\Gamma}} = f^*( \Omega^1_{C_{\Gamma}} \otimes 
\omega ).  
\end{equation}
The formula (\ref{eq:canonical}) combined with the Shimura isomorphism 
(\ref{eq:Shimura}) then give rise to the isomorphisms
\begin{equation}
H^0(C_{\Gamma}, \Omega^1_{C_{\Gamma}}) \simeq S_2(\Gamma) , 
\quad H^0(S_{\Gamma}, K_{S_{\Gamma}}) \simeq  
H^0(C_{\Gamma}, \Omega^1_{C_{\Gamma}} 
\otimes \omega) \simeq S_3(\Gamma).  
\end{equation}

Suppose that $S_{\Gamma}$ is a rational surface. 
Then $C_{\Gamma}$ is also a rational curve, so that
$\dim S_2(\Gamma) = 0$. We also have 
$\dim H^0(S_{\Gamma}, K_{S_{\Gamma}}) = 0$, and
hence $ \dim S_3(\Gamma) = 0$.  

Conversely, assume that $\dim S_k(\Gamma)=0 $ for $k =2, 3$.  
Then again the Shimura isomorphism (\ref{eq:Shimura}) implies that 
$H^0(C_{\Gamma}, \Omega^1_{C_{\Gamma}}) = 0$.
Hence $C_{\Gamma} \simeq \BP^1$.  Since $f$ is not isotrivial, the  
formula for the canonical bundle asserts that we may write 
$\omega$ as $\omega = \cO_{\BP^1}(a)$ and $K_{S_{\Gamma}} = f^*(\cO_{\BP^1}(a-2))$ 
with a suitable positive integer $a$. 
Note that the positivity of $a$ follows from the positivity of the direct image of the relative 
canonical sheaf  $\omega = f_*(\omega_{S_{\Gamma}/C_{\Gamma}})=\omega$.

Under this situation, $\dim S_3(\Gamma)= \dim H^0(S_{\Gamma}, K_{S_{\Gamma}}) = \dim H^0(\BP^1, 
\cO_{\BP^1}(a-2)) = 0 $ implies that $a = 1$.  
We also have the isomorphisms
$$
f_*\cO_{S_{\Gamma}} \simeq \cO_{\BP^1}, \quad 
R^1f_*\cO_{S_{\Gamma}} \simeq \omega^{\vee} \simeq \cO_{\BP^1}(-1). 
$$
Further the Leray spectral sequence shows that $H^1(S_{\Gamma}, \cO_{S_{\Gamma}})=0$.  
Moreover, for every $k > 0$, $ K^k_{S_{\Gamma}} = f^*(\cO_{\BP^1}(-k)) $ 
has no section. These facts finally imply that $S_{\Gamma}$ is a rational 
surface.  
{\it q.e.d}
\smallskip

Next, let us consider the self-fiber product of the elliptic modular surface 
associated to $\Gamma$
$$
f^2=f\times f: S_{\Gamma}^2=S_{\Gamma}\times_{C_{\Gamma}} S_{\Gamma}\lra C_{\Gamma}
$$
In general, $S_{\Gamma}^2$ has singularities which arise from the critical points of 
the map $f:S_{\Gamma} \lra C_{\Gamma}$.  However, if $f$ has only semistable fibers 
(that is, the singular fibers are all reduced and have only nodal singularities), 
then the fiber product $S^2_{\Gamma}$ has only ordinary double points as its
singularities.  
\smallskip

\begin{Proposition}\label{prop:canonical}
Suppose that $f:S_{\Gamma} \lra C_{\Gamma}$ has only semistable fibers. Then there exists a small projective resolution 
$\pi: \tilde{S}^2_{\Gamma}\lra S^2_{\Gamma}$. Moreover, if we set 
$h = f^2 \circ \pi$, then the canonical bundle of $\tilde{S}^2_{\Gamma}$
can be written  as 
\begin{equation}\label{eq:can}
K_{\tilde{S}^2_{\Gamma}}\simeq (h)^*(\Omega^1_{C_{\Gamma}} \otimes \omega^2).
\end{equation}
\end{Proposition} 
\smallskip

{\it Proof.}  The first assertion follows from Lemma 3.1 \cite{Sch2}.    Since 
$\pi$ is a small resolution, $K_{\tilde{S}^2_{\Gamma}}=\pi^*(K_{S^2_{\Gamma}})$.  
On the other hand, since $S^2_{\Gamma}$ is 
a hypersurface in $S_{\Gamma} \times S_{\Gamma}$ (i.e., it is the pull-back of 
the diagonal $\triangle : C_{\Gamma} \hookrightarrow C_{\Gamma} \times C_{\Gamma}$),
we can derive from the adjunction formula that 
\begin{eqnarray*}
K_{S^2_{\Gamma}} & = &  
(K_{S_{\Gamma} \times S_{\Gamma}} + 
f^{-1}(\triangle(C_{\Gamma})))_{|f^{-1}(\triangle(C_{\Gamma}))}  \\ 
 & = & (f^2)^* (( \Omega^1_{C_{\Gamma}})^{\otimes 2}
\otimes \omega^2 \otimes ( \Omega^1_{C_{\Gamma}})^{-1})  \\
 &  =  & (f^2)^* (\Omega^1_{\Gamma} \otimes 
\omega^2).
\end{eqnarray*}
This gives rise to the formula ({\ref{eq:can}}).  {\it q.e.d.}
\smallskip

\begin{Corollary} (\rm{cf. \cite{Sok}.})
 Under the assumption of Proposition \ref{prop:canonical}, we have a canonical 
isomorphism
$$
H^0(\tilde{S}^2_{\Gamma}, K_{\tilde{S}^2_{\Gamma}}) 
\simeq H^0(C_{\Gamma}, \Omega^1_{C_{\Gamma}} 
\otimes \omega^2 ) \simeq  S_4(\Gamma). 
$$
\end{Corollary}
\smallskip

For the full cohomology group $H^3(\tilde{S}^2_{\Gamma}, \Q)$, the 
following fact holds (cf. \S 1, \cite{Sch1}).  

\begin{Lemma}\label{lem:h3}
 {\rm (Lemma 1.7, \cite{Sch1}.)} 
The rational Hodge structure of 
$ H^3(\tilde{S}^2_{\Gamma},\Q)$ is 
isomorphic to the direct sum of three copies of 
$H^1(C_{\Gamma}, \Q)[-1]$ and a piece of pure 
type $(3, 0), (0,3)$.  In particular,  if 
$C_{\Gamma}$ is a rational curve, then 
$H^3(\tilde{S}^2_{\Gamma}, \Q)$ is of pure type 
$(3,0), (0,3)$.  
\end{Lemma}
\smallskip

\begin{Corollary}
Assume that $S_{\Gamma} \lra C_{\Gamma}$ is a rational elliptic modular 
surface associated to $\Gamma$ with only semistable singular fibers.  
Then any small resolution $\tilde{S}^2_{\Gamma}$ of the self-product 
$S^2_{\Gamma}$ is a rigid Calabi--Yau threefold. 
\end{Corollary}
\smallskip

{\it Proof.} 
Since $f:S_{\Gamma} \lra C_{\Gamma}$ is a rational surface, 
it is immediate that $C_{\Gamma} \simeq \BP^1$.  
Let $h:\tilde{S}^2_{\Gamma}\lra \BP^1$ be the natural fibration.  
We may write $\omega$ as $\omega = \cO_{\BP^1}(1)$ using the formula 
(\ref{eq:can}).  Then we see that 
$$
K_{\tilde{S}^2_{\Gamma}} = h^*(\Omega^1_{\BP^1} \otimes \omega^2)=
h^*(\cO_{\BP^1} (-2 +2)) = \cO_{\tilde{S}^2_{\Gamma}}. 
$$
Recall from Lemma (\ref{lem:h3}) that the rational Hodge structure 
$H^3(\tilde{S}^2_{\Gamma}, \Q) $ is of pure type 
$(3,0), (0,3)$.  So it follows that 
$H^1(\tilde{S}^2_{\Gamma}, \Omega^2_{\tilde{S}^2_{\Gamma}})=
H^1(\Theta_{\tilde{S}^2_{\Gamma}})=0.$  
We also see (cf. \cite{Sch1}) that 
$h^1(\cO_{\tilde{S}^2_{\Gamma}})=h^2(\cO_{\tilde{S}^2_{\Gamma}})=0$.
These facts then imply that $\tilde{S}^2_{\Gamma}$ is a rigid 
Calabi--Yau threefold.  
{\it q.e.d.}
\smallskip

\begin{Remark}{\rm There are in total six arithmetic subgroups $\Gamma 
\subset \text{SL}_2(\Z)$ which give rise to rational elliptic modular surfaces  
with only semistable fibers. The list of such groups coincides with the list of 
Beauville \cite{B}.  (Beauville classified elliptic surfaces with exactly 
four semistable fibers (i.e., of type $I_b,\,b\geq 1$) over $\BP^1$.)  
\smallskip

In Table \ref{tab:beauville}, groups $\Gamma$ are listed in the first column,
4 singular fibers of type $I_b,\,b\geq 1$ in the second column
and $h^{1,1}$ (resp. the Euler characteristic $\chi$) of the corresponding Calabi--Yau 
threefolds $\tilde{S}^2_{\Gamma}$ in the third (resp. fourth) column.
Note that the Euler characteristic $ \chi(\tilde{S}^2_{\Gamma})$ is exactly
equal to $2 h^{1,1}$. This is because $\tilde{S}^2_{\Gamma}$ is rigid.
Also note that for any group $\Gamma$ in Table \ref{tab:beauville}, 
$\dim S_4(\Gamma) = 1$. This follows from the isomorphism $S_4(\Gamma) \simeq 
H^{3,0}(\tilde{S^2_{\Gamma}}) \simeq \C  $.
\smallskip

\begin{table}[h]
\caption{}
\label{tab:beauville}
\begin{center}
\begin{tabular}{l|cccc|c|c} 
\quad \quad \quad  $\Gamma$ &    &  &  & &  $h^{1,1}(\tilde{S}^2_{\Gamma})$ & 
$\chi(\tilde{S}^2_{\Gamma})$ \\ \hline 
$\Gamma(3)$ & $I_3$ & $I_3$ & $I_3$ & $I_3$ & 36 & 72 \\
$\Gamma_1(4) \cap \Gamma(2)$ & $I_4$ & $I_4$ & $I_2$ & $I_2$ & 40 & 80  \\
$\Gamma_1(5)$  & $I_5$ & $I_5$ & $I_1$ & $I_1$  &52 & 104 \\
$\Gamma_1(6)$ & $I_6$ & $I_3$ & $I_2$ & $I_1$ &  50 & 100 \\
$\Gamma_0(8) \cap \Gamma_1(4)$ & $I_8$ & $I_2$ & $I_1$ & $I_1$  & 70 & 140 \\
$\Gamma_0(9) \cap \Gamma_1(3)$ & $I_9$ & $I_1$ & $I_1$ & $I_1$ & 84  & 168  
\end{tabular}
\end{center}
\end{table}
}
\end{Remark} 
\smallskip

Now we quote a theorem which is a special case of the fundamental 
result due to Sato--Kuga--Shimura (\cite{KS}; cf. \cite{D}). 
We make use the argument of \cite{D}, where 
Deligne worked over the case when $\Gamma = \Gamma(N)$ (full modular case).  However his argument works 
for other cases like $\Gamma = \Gamma_1(N)$
as far as the universal family of generalized 
elliptic curves exists. 

\begin{Theorem} \label{thm:self}
Let $\Gamma$ be an arithmetic subgroup of $\text{SL}_2(\Z)$ listed in Table 
\ref{tab:beauville}, and let  
$f:S_{\Gamma} \lra \BP^1$ the rational elliptic
surface associated to $\Gamma$.  Let $f^2:\tilde{S}^2_{\Gamma}\lra \BP^1$  
be the smooth rigid Calabi--Yau threefold arising from the self-product
$S^2_{\Gamma}$. Then $\tilde{S}^2_{\Gamma}$ is modular,  
that is, the L-series of $\tilde{S}^2_{\Gamma}$ coincides 
with the Mellin transform of a weight $4$ newform of $\Gamma$ up to finite Euler factors coming from the bad primes:
$$
L(\tilde{S}^2_{\Gamma}, s) = L(f, s) \quad \mbox{for} 
\quad f \in S_4(\Gamma).   
$$ 
\end{Theorem}
\begin{Remark}\label{rem:bad}
The bad primes for $\tilde{S}^2_{\Gamma}$ depends on the level of the discrete group $\Gamma$.  For example, we can easily see that 
only $2$ and  $3$ 
are the bad primes for $\ \Gamma_1(6)$.  
\end{Remark}
\medskip

\section{\bf The geometry of Verrill's Calabi--Yau threefold}

\hskip 1cm In this section, we will look into the geometric structure of Verrill's 
rigid Calabi--Yau threefold $X$ constructed in Section 3.
Note that from now on all varieties are defined over $\Q$. 

We denote by $[x:y:z:w]$ the homogeneous coordinate of $\BP^3$ and by $t$ the 
inhomogeneous coordinate of $\BP^1$.  Then Verrill's rigid Calabi--Yau threefold 
$X$ is a desingularization of the hypersurface in $\BP^3 \times \BP^1$:

\begin{equation}\label{eq:verrill}
X' = \{ (x+y+z+w)(x^{-1}+y^{-1}+z^{-1}+w^{-1}) - \frac{(t-1)^2}{t}+4 = 0 \} 
\subset \BP^3\times \BP^1.
\end{equation}
It is easily seen that this equation is equivalent to the following equation:
\begin{eqnarray}\label{eq:v2}
F(x, y, z, w, t)  & = &  (x + y + z + w) 
\cdot [(y z+ z x + x y) w + xyz] \cdot t \\
 &  & \quad \quad  -  (t +1)^2 \cdot x y z w = 0. \nonumber
\end{eqnarray}

\noindent{\bf The elliptic modular surface associated to $\Gamma_1(6)$.}

\vspace{0.5cm}
Let us consider the hypersurface $S$ in $\BP^1 \times \BP^2$ 
defined by
\begin{eqnarray}\label{eq:ell}
H(x, y, z, s) = (s+1)xyz - (x+y+z)(yz+zx+xy) = 0.
\end{eqnarray}
The hypersurface $S$ has three singular points at 
$$
[s, [x:y:z]]= [\infty, [1:0:0]], \ [\infty, [0:1:0]],  \  [\infty, [0:0:1]] .
$$ 
Let 
\begin{equation*}
 \pi: \tilde{S} \lra S
\end{equation*}
be the blowing up of these three points.  
 Then it is easy to see that  
$\tilde{S}$  is the minimal resolution  of $S$ 
and $f= p_1 \circ \pi: \tilde{S} \lra \BP^1$ induces the 
structure of the elliptic surface.  
Moreover the rational elliptic surface 
 $f:\tilde{S} \lra \BP^1$ 
has only four semistable singular fibers 
of types $I_1, I_2, I_3, I_6$ (see Table \ref{tab:sing}):   
\begin{table}[h]
\caption{}\label{tab:sing}
\begin{center}
\begin{tabular}{|c||c|c|c|c|}\hline
$s$ & $8$  & $-1$ & $0$ & $\infty$ \\ \hline 
singular fiber 
& $I_1$ & $I_2$ & $I_3$ 
& $I_6$ \\   \hline 
\end{tabular}
\end{center}
\end{table}

By \cite{B}, $ f:\tilde{S} \lra \BP^1 $ is an integral model 
of the elliptic modular surface 
$$ 
f:S_{\Gamma_1(6)} \lra \overline{H/\Gamma_1(6)} \simeq \BP^1
$$   
of $\Gamma_1(6)$. 
\smallskip

\begin{Theorem}\label{thm:main}
Verrill's Calabi--Yau threefold $X$ defined in 
(\ref{eq:verrill}) is birationally equivalent over $\Q$ 
to a crepant resolution $\tilde{S}^2_{\Gamma_{1}(6)}$ of the self-product of 
the elliptic modular surface $S_{\Gamma_{1}(6)}$.  More precisely, there exists 
a birational map defined over $\Q$ between $X$ and $\tilde{S}^2_{\Gamma_{1}(6)}$. 
\end{Theorem}
\smallskip

\begin{Remark}{\rm There are two fibrations, one is the fibration of $X \lra \BP^1$ 
of K3 surfaces and the other is the fibration 
$f^2: \tilde{S}^2_{\Gamma_1(6)}\lra \BP^1$ of abelian 
surfaces.  We will show below that these two fibrations are not equal. 
(However, we do not know any relation (apart from being non-equal) 
between these two fibrations.)}
\end{Remark}
\medskip

\noindent{\bf Proof of Theorem \ref{thm:main}:  Birational transformations.}

To prove Theorem \ref{thm:main}, we will construct a birational 
map over $\Q$ between the two varieties explicitly.  
We define the birational map: 

$$
\begin{array}{ccc}
\pi: \BP^2 \times \BP^2 &  \longrightarrow &  
\BP^3 \times \BP^1 \\ 
  &  &  \\
  \ [  x : y  : z  ]  \  \times  \ [T : W : U]  & 
    \mapsto & [ x : y : z : w ] \times (t)
\end{array} 
$$ 

\noindent by putting 

$$
x = x, \  y = y, \  z = z,  \  t = T/U,  \ w = W \cdot (x + y +z)/U. 
$$

\noindent Pulling back the equation 
$F(x, y, z, w, t)$ in (\ref{eq:v2}) via $\pi$, we obtain 
the equation 
$$
\begin{array}{rcl}
\tilde{F}:  & = &   (x+y+z)^2(yz + xz + xy)TW^2  \\
     &  & - (x + y + z) \cdot[(x + y +z)(yz + xz + xy) - x y z ] T W U  \\ 
 & &  + x y z (x + y + z) T U^2 -(x + y + z) x  y  z 
 [T^2 + U^2]  W = 0
\end{array}
$$
Consider the hypersurface $X'' \subset \BP^2 \times \BP^2$  defined 
by $\tilde{F}  = 0$. Of course, $X''$ is birationally equivalent to 
the original hypersurface $X'$ in (\ref{eq:verrill}).  
We have the following commutative diagram:
$$
\begin{array}{ccc}
X''          & \stackrel{\iota}{\hookrightarrow} & \BP^2 \times \BP^2 \\
   &  &  \\
f_1 \downarrow \quad & \swarrow p_1      &        \\
  &  &  \\
\BP^2          &     &      \\
\end{array}
$$

Now recall the hypersurface $ S \subset \BP^1 \times \BP^2 $ defined in (\ref{eq:ell}). 
Considering  the birational map $ \mu= 
(p_{2})_{|S}: S \lra \BP^2 $  
and taking  the pull-back morphism $ f_1 $ by $\mu$, we obtain 
the following commutative diagram:
$$
\begin{array}{ccccc}
\BP^1 \times \BP^2 \times \BP^2 & \supset
 &  X'''   & \stackrel{\nu}{\lra} &  X'' \\
& &     &   &    \\
\downarrow p_{12} & & f'_1  \downarrow \quad  &  & f_1 \downarrow \quad \\
   &  &   &   \\
\BP^1 \times \BP^2 &  \supset  & S  & \stackrel{\mu}{\lra}  &  \BP^2  
\end{array} 
$$
Here we set $X''' := S \times_{\BP^2} X''$.
Note that the natural induced map $\nu:X''' \lra X''$ is a 
birational morphism.  Now we see that $X'''$ is a complete intersection 
  defined by the equations 
$$
\begin{array}{ccl}
\tilde{F} & = & 0  \\
H(x, y, z) & = & (s+1) (x y z) -  (x + y + z)(yz + zx + xy) = 0.
\end{array}
$$
Now using the second equation,  $\tilde{F}= 0 $ can  
be transformed into the following equation:
$$
\begin{array}{rcl}
\tilde{F}:  & = &   (x+y+z) \cdot (s+1)xyz\cdot T W^2  \\
     &  & + (x + y + z) \cdot[(s+1)xyz - x y z ] T 
       W U  \\ 
 & &  + x y z (x + y + z) T U^2 -(x + y + z)[T^2 + U^2] x  y  z  W  \\
      & = & (x+y+z) xyz \times \{ (s+1) \cdot T W^2 -s  T W U + T U^2 - T^2 W -W U^2  \}  
\end{array}
$$

\noindent Forgetting the factor $(x+y+z) xyz $ from $\tilde{F}$, 
we obtain the equation  
$$
G(T, W, U) =  (s+1) \cdot T W^2 + s  T W U + T U^2 - T^2 W -W U^2   = 0.
$$
Therefore we see that the subvariety $X'''$ is 
isomorphic to the complete 
intersection in $\BP^1 \times \BP^2 \times \BP^2$  defined by 
the equations
$$
\begin{array}{ccl}
G(T, W, U) &  = &   (s+1) T W^2 + s  T W U + T U^2 - T^2 W -W U^2 = 0, \\
H(x, y, z) &  = &  (s+1) x y z - (x + y + z)(yz + zx + xy) = 0.
\end{array}
$$

Next, considering  the  the birational transformation of $\BP^2$ given by 
$$
T= Z \cdot (X+Y+Z), \quad W = -X Y, \quad U = Y(X+Y+Z).  
$$  

\noindent Then the polynomial  $G(T, W, U) $ above can be transformed into the polynomial

$$
\tilde{G}(X, Y, Z) =( s+1) XYZ -(X+Y+Z)(YZ+ZX+XY), 
$$  

\noindent which is nothing but $H(X, Y, Z) $. 
Hence the complete intersection 
$ X''' \subset \BP^1 \times \BP^2 \times \BP^2$ 
is  isomorphic to the fiber product 
$S \times_{\BP^1} S \lra \BP^1$ 
of the elliptic surface $ f:S \lra \BP^1$ defined by
the equation $H(x, y, z) = 0$.  Hence, $X'$ defined in (\ref{eq:verrill}) 
is birationally equivalent to the self-product $S \times_{\BP^1} S$ of
$S$. Since the minimal resolution of $S$ is isomorphic to the elliptic 
modular surface $S_{\Gamma_1(6)}$, $X$ is birationally
equivalent to the self product $\tilde{S}^2_{\Gamma_1(6)}$. 
Observing that the birational maps discussed above are all defined over $\Q$, 
this completes the proof of Theorem \ref{thm:main}. 
{\it q.e.d.}
\medskip

\section{\bf Proof for the modularity of $X$ (Theorem \ref{thm:modular}) }

Since the birational map between Verrill's example $X$ and $\tilde{S}^2_{\Gamma_1(6)}$ 
constructed in Section 5 is defined over $\Q$, we obtain the following Lemma.

\begin{Lemma}\label{lem:compat}
The birational map constructed in Section 5 gives rise to the following 
isomorphism compatible with the Galois action up to  bad primes $2$ and  $3$
$$
 H^3_{et}(\bar X,\Ql)\simeq H^3_{et}(\bar{S}^2_{\Gamma_1(6)}, \Ql).
$$
Consequently, we have
$$L(X,s) = L(\tilde{S}^2_{\Gamma_1(6)}, s).$$
up to Euler factors at bad primes $2$ and $3$.  
\end{Lemma}
\smallskip

For a subgroup $\Gamma\subset\text{SL}_2(\Z)$, its projectivization
$\Gamma/\{\pm I_2\}$ will be denoted by $P\Gamma$.
\medskip
By virtue of Theorem \ref{thm:self}, Theorem \ref{thm:main}, Remark \ref{rem:bad} and 
Lemma \ref{lem:compat}, the proof of 
Theorem \ref{thm:modular}  (the modularity of Verrill's Calabi--Yau threefold $X$)  
is reduced to the proof of the following lemma.

\begin{Lemma} 
The projective image of $\Gamma_1(N)$ and $\Gamma_0(N)$ are the same
in $\text{PSL}_2(\Z)$ if and only if $N$ is a divisor of $4$
or a divisor of $6$.  In particular, ${\BP}\Gamma_1(6)={\BP}\Gamma_0(6)$
and so it follows that
$S_4({\BP}\Gamma_1(6))\simeq S_4({\BP}\Gamma_0(6))$.
\end{Lemma}
\smallskip

{\it Proof}. Let $\pmatrix a & b \\ c & d \endpmatrix$ be a matrix in
$\Gamma_0(N)$ or $\Gamma_1(N)$.  Then $ad\equiv 1\pmod N$.  
The first assertion follows from the fact that the only
positive integers $N$ for which $ad\equiv 1\pmod N$ implies that $a\equiv d
\equiv \pm 1\pmod N$ are the divisors of $4$ and $6$. {\it q.e.d.} 
\smallskip

\begin{Remark}{\rm Verrill \cite{V} has an alternative proof of the modularity for $X$.
Her proof is along the lines of Livn\'e's paper \cite{Li1}.  It makes use of
the Serre criterion \cite{Sr} and Faltings results \cite{F} to prove equality
of two $L$--series.  The main point was to show that {\it finitely many} 
Euler factors of two $L$--series coincide.}
\end{Remark}
\smallskip

\begin{Remark}
{\rm Our geometric proof of the modularity conjecture for $X$ 
guarantees that the $2$--dimensional Galois representations 
associated to $X$ and that attached to the modular form are indeed isomorphic (i.e.,
have the same semi-simplification).  This provides a confirmation
to the conjecture of Fontaine and Mazur \cite{FM}.}
\end{Remark}
\smallskip

We close this section by posing an open problem.

\begin{Problem}{\rm 
May one use the method of Wiles to establish the modularity of $X$?  
More concretely, find a single good prime $\ell$ and establish the modularity for
the residual mod $\ell$ Galois representation associated to the rigid
Calabi--Yau threefold $X$ in question.}  
\end{Problem}
\medskip

\section{\bf The intermediate Jacobians of rigid Calabi--Yau threefolds}

\hskip 1cm In this section, we shall define the intermediate Jacobian for 
Calabi--Yau threefolds, following the exposition of Bloch [Bl].  
\medskip

\begin{Definition}{\rm 
Let $X$ be a smooth projective Calabi--Yau
threefold defined over $\C$.  There is a Hodge filtration 
$$H^3(X,\C)=F^0\supset F^1\supset F^2\supset F^3\supset (0)$$
where 
$$F^1=H^{3,0}\oplus H^{2,1}\oplus H^{1,2}$$
$$F^2=H^{3,0}\oplus H^{2,1}$$
$$F^3=H^{3,0}.$$
This complex vector space also has a $\Z$--structure defined by
the image $H^3(X,\Z)\to H^3(X,\C)$.
\smallskip

Now we recall some key properties of Hodge filtrations not only
for $H^3$ but for more general cases. We write $\cF$ 
for generic filtrations.  We have
$$\cF^i\oplus \overline{\cF^{2r-i}}\simeq H^{2r-1}(X,\C),\quad
\cF^i\cap\overline{\cF^{2r-i-1}}\simeq H^{i,2r-i-1}$$
where $H^{i,j}\simeq H^j(X,\Omega_X^i)$ and $r\in\{1,2,3\}$.  
{F}rom this, one can compute that
$$\text{dim}_{\C} \cF^r=\frac{1}{2}\text{dim}_{\C} H^{2r-1}(X,\C)$$
and
$$\cF^r\cap\text{Image}(H^{2r-1}(X,\Z)\to H^{2r-1}(X,\C))=(0).$$
Therefore, the quotient
$$J^r(X)=H^{2r-1}(X,\C)/(\cF^r + H^{2r-1}(X,\Z))$$
is a compact complex torus, called the {\it intermediate Jacobian}.
By the Poincar\'e duality, we have
$$H^{2r-1}(X,\C)/\cF^r\simeq \cF^{4-r}\,H^{7-2r}(X,\C)^*$$
where $*$ denotes $\C$--linear dual.  It then follows that
$$J^r(X)\simeq \cF^{4-r}\,H^{7-2r}(X,\C)^*/H_{7-2r}(X,\Z).$$}
\end{Definition}
\smallskip

{\bf { Digression}.}  
Now assume that $X$ is rigid, and take $r=2$ and $\cF=F$. Then 
we have
$$F^1=H^{3,0}=F^2=F^3$$
and
$$J^2(X)\simeq H^3(X,\C)/(F^2+H^3(X,\Z))
\simeq H^{3,0}(X)^{*}/H_3(X,\Z)$$
is a complex torus of dimension one.  This means that there is an
elliptic curve $E$ such that $E(\C)\simeq J^2(X)$.
We will formulate the following question.
\smallskip

\begin{Question}
Is it true that a rigid Calabi--Yau threefold
defined over $\Q$ is modular if and only if the intermediate
Jacobian is an elliptic curve defined over $\Q$?
\end{Question} 
\vfill

\pagebreak

{\bf Acknowledgments}

We thank several people who read through the earlier versions
of this paper, and gave us constructive comments and suggestions.
These include I. Chen, K. Khuri-Makdisi, R. Livn\'e,
A. Sebbar, and J.-P. Serre.

At the early stages of this work, M.-H. Saito was a visiting
Professor at Queen's University in the summer of 1998.  
The visit was made possible by the grant from JSPS.  

During the course of this work, Yui was a visiting Professor 
at various universities and research institutes, e.g., 
CRM Montr\'eal, The Fields Institute, MSRI Berkeley, 
Max-Planck-Institut Bonn, Chiba University, Tsuda College, 
Kyoto University and T\^ohoku University. 
She is grateful for hospitality and supports of all these
institutions.

\end{document}